


\input amstex
\documentstyle{amsppt}

\magnification=1200

\overfullrule=0pt

\hsize=165truemm

\vsize=227truemm


\def\ga#1{{\accent"12 #1}}

\def\p#1#2{{{\Bbb P}^{#1}_{#2}}}

\def\Pic{\operatorname{Pic}}

\def\im{\operatorname{im}}

\def\rk{\operatorname{rk}}

\def\pd{\operatorname{pd}}

\def\ann{\operatorname{ann}}

\def\ext{\operatorname{Ext}}

\def\coker{\operatorname{coker}}

\def\depth{\operatorname{depth}}

\def\hom{\operatorname{Hom}}

\def\Hom#1{{{\Cal H}\kern -0.25ex{\italic om\/}_{\Ofa#1}}}

\def\Syz{{\Cal S}\kern -0.25ex{\italic yz\/}}

\def\Ofa#1{{{\Cal O}_{#1}}}

\def\Sim{\mathop{{\Cal S}\/}\nolimits}

\def\Ext#1#2{{{\Cal E}\kern -0.25ex{\italic
xt\/}^#1_{\Ofa{#2}}}}

\def\mapright#1{\mathbin{\smash{\mathop{\longrightarrow}
\limits^{#1}}}}

\def\mapdown#1{\Big\downarrow\rlap{$\vcenter{\hbox
{$\scriptstyle#1$}}$}}


\topmatter

\title
Even sets of nodes are bundle symmetric
\endtitle

\author
G. Casnati and F. Catanese
\endauthor

\address
Gianfranco Casnati: Dipartimento di Matematica Pura ed
Applicata, Universit\ga a degli Studi di Padova,
via Belzoni 7, I--35131 Padova (Italy)
\endaddress

\email
casnati\@pdmat1.math.unipd.it
\endemail

\address
Fabrizio Catanese: Dipartimento di Matematica, Universit\ga a
degli Studi di Pisa, via Buo\-nar\-ro\-ti 2, I--56127 Pisa
(Italy)
\endaddress

\email
catanese\@dm.unipi.it
\endemail

\endtopmatter

\document

\head
0. Introduction
\endhead

Let $k$ be an algebraically closed field of characteristic
$p\ne2$ and let $F:=\lbrace f=0\rbrace\subseteq\p3k$ be a
normal surface of degree
$d$. Let
$\pi\colon \widetilde{F}\to F$ be a minimal resolution of
singularities. We denote by $H\subseteq F$ a general plane
section of $F$ defined by a general linear form $h$. Assume,
for simplicity, that $F$ is a nodal surface (i.e. its
singularities are only ordinary quadratic, nodes for short).

Let $\Delta$ be a subset of the set of nodes of $F$ and let
$\widetilde{\Delta}:=\pi^{-1}(\Delta)$. $\Delta$ is said to
be a {\sl
$\delta/2$--even set of nodes}\/, $\delta=0,1$, if the class
of
$\widetilde{\Delta}+\delta \pi^*H$ in $\Pic(\widetilde{F})$
is
$2$--divisible (when $\delta=0$ we shall simply say that
$\Delta$ is even).

If $\Delta$ is $\delta/2$--even then there exists a double
cover
$\widetilde{p}\colon \widetilde{S}\to\widetilde{F}$ branched
exactly along $\widetilde{\Delta}+\delta \pi^*H$ and (cf.
[Ct] 2.11, 2.13) it is possible to blow down
$\widetilde{p}^{-1}(\widetilde{\Delta}+\delta \pi^*{H})$
getting a commutative diagram
$$
\matrix
\widetilde{S}&\mapright{\widetilde{\pi}}&S\\
\mapdown{\widetilde{p}}&&\mapdown p\\
\widetilde{F}&\mapright{\pi}&F
\endmatrix\tag0.1
$$
where $S$ is a nodal surface and $p$ is finite of degree $2$
branched exactly on $\Delta$ when $\delta=0$ respectively
$\Delta$ and $H$ when $\delta=1$ (in this case $\delta$ has
to be even). The surface
$S$ is then endowed with a natural involution $i$ such that
$F\cong S/i$ and $p$ is the quotient map. Thus we have an
$\Ofa F$--linear map $i^\#\colon p_*\Ofa S\to p_*\Ofa S$
giving rise to a splitting of $\Ofa F$--modules
$p_*\Ofa S\cong\Ofa F\oplus{\Cal F}$ where $\Ofa F$ and
$\Cal F$ are the $+1$ and $-1$ eigenspaces of $i^\#$. Notice
that
${\Cal F}_{\vert F\setminus\Delta}$ is invertible, hence
reflexive.

The multiplication map $\Ofa S\times\Ofa S\to\Ofa S$
induces a symmetric non--degenerate bilinear form ${\Cal
F}\times{\Cal F}\to\Ofa F$. If $\delta=1$ such a map factors
through the multiplication by $h$, $\lambda(h)\colon\Ofa
F(-1)\to\Ofa F$. Thus, in both cases, we get a symmetric
bilinear map ${\Cal F}\times{\Cal F}\to\Ofa F(-\delta)$
inducing a monomorphism
$\sigma\colon {\Cal
F}(\delta)\rightarrowtail\check{\Cal F}:=\Hom{F}\big({\Cal
F},\Ofa F\big)$ which is obviously an isomorphism outside
$\Delta$, whence a global isomorphism since $S$ is normal.

Thus  $\Cal F$ is isomorphic to ${\Cal F}\check{\ }\check{\
}$ and, indeed, $\Cal F$ is reflexive  which means that the
natural map  ${\Cal F}\mapright\mu{\Cal F}\check{\ }\check{\
}$ is an isomorphism (infact
$\coker(\mu)=0$ since $\Cal F$ and ${\Cal F}\check{\
}\check{\ }$ have the same Hilbert polynomial).

If $x\in\Delta$ the completion $\widehat{\Cal
F}_x$ fits into the exact sequence
$$
0\to\widehat{\Cal O}_{\p3k,x}^{\oplus2
}\mapright{\mu}\widehat{\Cal
O}_{\p3k,x}^{\oplus2}\to\widehat{\Cal F}_x\to0
$$
where $\mu$ is induced by a matrix of the form
$$
\pmatrix
w&z\\
z&y
\endpmatrix,
$$
$w,z,y$ being local parameters in $\widehat{\Cal
O}_{\p3k,x}$. Therefore we see that
${\Cal F}_x$ is Cohen--Macaulay as a module over
$\Ofa{\p3k,x}$.

We conclude that $\Cal F$ is a $\delta/2$--quadratic sheaf
in the sense of the following definition.

\definition{Definition 0.2}
Let $X$ be a locally Cohen--Macaulay scheme. We say that a
reflexive, coherent, locally Cohen--Macaulay
$\Ofa X$--sheaf $\Cal F$ is a
$\delta/2$--quadratic sheaf on $X$, $\delta\in\Bbb Z$, if
there exists a symmetric isomorphism $\sigma\colon{\Cal
F}(\delta)\mapright\sim\Hom{X}\big({\Cal F},\Ofa X\big)$.
\enddefinition

The aim of sections 1 and 2 is to prove in dimension $3$ the
following characterization of quadratic sheaves on
hypersurfaces in projective space.

\proclaim{Theorem 0.3}
Let $F\subseteq\p3k$ be a  surface of degree $d$ and let
$\Cal F$ be a
$\delta/2$--quadratic sheaf on $F$. Then
$\Cal F$ fits into a sequence of the form
$$
0\to\check{\Cal E}(-d-\delta)\mapright\varphi{\Cal E}\to{\Cal
F}\to0\tag0.3.1
$$
where $\Cal E$ is a locally free $\Ofa{\p3k}$--sheaf and
$\varphi$ is a symmetric map.
\qed
\endproclaim

An entirely analogous proof with more complicated notations
gives the same result in all dimensions.

Theorem 0.3 and the above discussion
yield the following

\proclaim{Corollary 0.4}
Let $F\subseteq\p3k$ be a nodal surface of degree $d$. Then
every
$\delta/2$--even set of nodes $\Delta$ on $F$, $\delta=0,1$,
is the degeneracy locus of a symmetric map of locally free
$\Ofa{\p3k}$--sheaves
$\check{\Cal E}(-d-\delta)\mapright\varphi{\Cal E}$ (i.e.
$F$ is the locus where $\rk(\varphi)\le\rk{\Cal E}-1$,
$\Delta$ is the locus where
$\rk(\varphi)\le\rk{\Cal E}-2$).
\qed
\endproclaim

In the above setting we say that $\Delta$ is a {\sl
bundle--symmetric}\/ set of nodes. If it is possible to find
such an
$\Cal E$ which is the direct sum of invertible
$\Ofa{\p3k}$--sheaves, then we say that
$\Delta$ is a {\sl symmetric}\/ set of nodes (see [Ca]).

Corollary 0.4 was conjectured in 1979 independently by W.
Barth and the second author. Barth proved in [Ba1] that
bundle--symmetric sets are even while in [Ct] the converse
result was proved under a cohomological assumption which
gives a complete characterization of symmetric sets of nodes.

As soon as Walter's beautiful solution of Okonek's
conjecture came out, it was immediately clear that his method
would also work in our case.

Finally, in section 3 we apply corollary 0.4 to the study of
nodal surfaces of low degree $d$, namely $d=4,5$. This study
ties up to an interesting history for which we defer the
reader e.g. to [Gal], [Ct], [Bea], [Chm], [Mi], [Ba2],
[J--R].

\subhead
Acknowledgements
\endsubhead
We would like to thank C. Walter for sending his inspiring
preprint [Wa] to the second author in June 1994.

Both the authors acknowledge support from the AGE project
H.C.M. contract ERBCHRXCT 940557 and from 40 \% M.U.R.S.T..

The final version was written while the second author was
\lq\lq Professore distaccato\rq\rq at the Accademia dei
Lincei.

\head
1. A resolution of $\Cal F$
\endhead

In this section we deal with the construction of a
resolution ${\Cal R}^*$ of any $\delta/2$--quadratic sheaf
$\Cal F$ on a  surface
$F\subseteq\p3k$ of degree $d$. We shall make use of
notations and results proved in [Wa] about Horrocks'
correspondence.

\proclaim{Lemma 1.1}
Let $F\subseteq\p3k$ be a surface and let ${\Cal F}$ be a
$\delta/2$--quadratic sheaf on $F$. Then
$\pd_{\Ofa{\p3k,x}}{\Cal F}_x=1$ for each $x\in F$ and
$\Ext{1}{F}\big({\Cal F},\Ofa F\big)=0$.
\endproclaim
\demo{Proof}
By definition $\depth{\Cal F}_x=\dim{\Cal F}_x=2$. Then the
equality
$\pd_{\Ofa{\p3k,x}}{\Cal F}_x=1$ follows from
Auslander--Buchsbaum formula taking account that the depths
of ${\Cal F}_x$ as $\Ofa {F,x}$--module and as
$\Ofa{\p3k,x}$--module coincide.
$\Ext{1}{F}\big({\Cal F},\Ofa F\big)=0$ follows from [H--K],
theorem 6.1.
\qed
\enddemo

Let $d+\delta$ be even. From the spectral
sequence of the $\roman{Ext}$'s, lemma 1.1 and Serre's
duality it follows the existence of isomorphisms for $i=1,2$
$$
\aligned
&H^i\big(F,{\Cal F}(m)\big)\cong\ext_{\Ofa F}^{2-i}\big({\Cal
F}(m),\omega_{F\vert k}\big)\check{}\cong\\
&H^{2-i}\big(F,\Hom{F}\big({\Cal F},\Ofa F(d-4-m)\big)\big)\check{}
\cong H^{2-i}\big(F,{\Cal F}(d-4-m+\delta)\big)\check{}.
\endaligned\tag 1.2
$$
In particular, taking $m:=(d-4+\delta)/2$, there exists a
non--degenerate alternating form
$$
\Phi\colon H^1\big(F,{\Cal F}((d-4+\delta)/2)\big)\times
H^1\big(F,{\Cal F}((d-4+\delta)/2)\big)\to H^2\big(F,\Ofa
F\big)
$$
and we denote by $U$ a fixed maximal isotropic subspace with
respect to
$\Phi$.

Define
$$
W:=\cases
\bigoplus_{m<(d-4+\delta)/2}H^1\big(F,{\Cal F}(m)\big)
&\text{ if
$d+\delta$ is odd},\\
\bigoplus_{m<(d-4+\delta)/2}H^1\big(F,{\Cal F}(m)\big)\oplus
U &\text{ if $d+\delta$ is even}.
\endcases
$$

As usual $H_*^i$ is the Serre functors associating to a
quasi--coherent sheaf $\Cal G$ the graded module
$H_*^i\big(\p3k,{\Cal G}\big):=\bigoplus_{n\in{\Bbb
Z}}H^i\big(\p3k,{\Cal G}(n)\big)$. Let $\Gamma_*:=H_*^0$ and
let $D^*:={\bold R}\Gamma_*$ be its right derived functor in
the derived category with differentials $\delta^i\colon
D^i\to D^{i+1}$ (then
$H_*^i\big(\p3k,{\Cal G}\big)$ is the $i$--th cohomology
module of the complex ${\bold R}\Gamma_*({\Cal G})$).

As in [Wa], section 2 one considers the truncation. Let
$D^*$ be as above, let $r,s\in\Bbb Z$ and let $W\subseteq
H^r\big(D^*\big)$ be a subspace. Then $W$ may be pulled back
to $\overline{W}$ satisfying
$\im(\delta^{s-1})\subseteq\overline{W}
\subseteq\ker(\delta^s)$. We then denote by
$\tau_{>r}\tau_{\le s,W}(D^*)$ (if $W=H^r\big(D^*\big)$ we
will omit it in the subscripts) the complex $C^*$ defined as
follows:
$$
C^i:=\cases
0&\text{ if $i\le r-1$ or $i\ge s+1$,}\\
D^i&\text{ if $r+1\le i\le s-1$,}\\
D^r/\ker(\delta^r)&\text{ if $i=r$,}\\
\overline{W}&\text{ if $i=s$.}
\endcases
$$

Let us then consider the truncated complex
$C^*:=\tau_{>0}\tau_{\le1,W}{\bold R}\Gamma_*({\Cal F})$. By
definition we have
$$
H^i\big(C^*\big)=\cases W &\text{if $i=1$},\\
0 &\text{elsewhere},
\endcases
$$
and there is a natural map
$$
\beta\colon C^*\to\tau_{>0}\tau_{\le2}{\bold R}\Gamma_*({\Cal
F}).
$$
Since $H^1_*\big(F,{\Cal F}\big)$, hence $W$, has finite
length we can apply to the complex $C^*$ the syzygy bundle
functor ([Wa], theorem 0.4. Cf. also the construction given
after corollary 2.8). We obtain a locally free sheaf
$\Syz(C^*)$ and a morphism of quasi--coherent sheaves
$$
\widetilde{\beta}\colon \Syz(C^*)\to \Cal F
$$
such that $\beta=\tau_{>0}\tau_{\le2}{\bold
R}\Gamma_*(\widetilde{\beta})$ ([Wa], proposition 2.10).

Let $Q:=\coker\big(H^0_*\big(\widetilde{\beta}\big)\big)$.
This means that we have an exact sequence of the form
$$
H^0_*\big(\p3k,\Syz(C^*)\big)
\mapright{H^0_*\big(\widetilde{\beta}\big)} H^0_*\big(F,\Cal
F\big)\to Q\to0.
$$
Let $d_1,\dots,d_r$ be the degrees of a minimal set of
generators of
$Q$. These generators lift to $H^0_*\big(F,\Cal F\big)$,
allowing us to define an epimorphism
$$
\gamma\colon {\Cal E}:=\Syz(C^*)\oplus\bigoplus_{i=1}^r
\Ofa{\p3k}(-d_i)\twoheadrightarrow{\Cal F}
$$
which is surjective on global sections.
By construction $\Cal E$ is locally free.

If ${\Cal K}:=\ker(\gamma)$ then we have an exact sequence
$$
{\Cal R}^*:\quad 0\to{\Cal K}\mapright d{\Cal E}
\mapright\gamma{\Cal F}\to0.\tag1.3
$$

\proclaim{Proposition 1.4}
In the sequence (1.3) above $\Cal K$ is locally free and
$\rk{\Cal K}=\rk{\Cal E}$.
\endproclaim
\demo{Proof}
We know from lemma 1.1 that for each $x\in\p3k$ one has an
exact sequence of the form
$$
0\to{\Cal K}'_x\to{\Cal E}'_x\to{\Cal F}_x\to0
$$
where ${\Cal K}'_x$ and ${\Cal E}'_x$ are free (and depend
upon
$x\in\p3k$). Moreover
$\ann_{\Ofa{\p3k,x}}{\Cal F}_x\ne0$ since $\Cal F$ is
supported on $F$ then
$\rk{\Cal K}'_x=\rk{\Cal E}'_x$ (see [Kp], theorem 195). The
statement now follows from [Kp], theorem A, chapter 4
(Schanuel's lemma).
\qed
\enddemo

Proposition 2.10 of [Wa] implies for $i=1,2$ that
$H^i\big(C^*\big)=H^i_*\big(\p3k,\Syz(C^*)\big)=
H^i_*\big(\p3k,{\Cal E}\big)$ and $H^i_*\big(\gamma\big)$
coincides with
$H^i_*\big(\widetilde{\beta}\big)$, thus the above
construction yields:

\roster
\item"i)" $H^0_*\big(\gamma\big)\colon H^0_*\big(\p3k,{\Cal
E}\big)\to H^0_*\big(F,{\Cal F}\big)$ is surjective by
construction;

\item"ii)" $H^1_*\big(\gamma\big)\colon H^1_*\big(\p3k,{\Cal
E}\big)\to H^1_*\big(F,{\Cal F}\big)$ is injective since
$H^1_*\big(\p3k,{\Cal E}\big)=H^1\big(C^*\big)=W$ (in
particular $H^1(\p3k,{\Cal E}((d-4+\delta)/2)\big)=U$);

\item"iii)" $H^2_*\big(\gamma\big)\colon
H^2_*\big(\p3k,{\Cal E}\big)\to H^2_*\big(F,{\Cal F}\big)$ is
zero since
$H^2_*\big(\p3k,{\Cal E}\big)=0$.
\endroster

{}From the above remarks taking the cohomology of the sequence
(1.3) we then get:

\roster
\item"iv)" $H^1_*\big(\p3k,{\Cal K}\big)=0$;

\item"v)" $H^2_*\big(\p3k,{\Cal K}\big)\cong
H^1_*\big(F,{\Cal F}\big)/W$.
\endroster

Recall that ${\Cal F}(\delta)\cong\Hom{F}\big({\Cal F},\Ofa
F\big)$. On the other hand one has an exact sequence
$$
{\Cal C}^*:\quad 0\to\Ofa{\p3k}(-d-\delta)
\mapright{\lambda(f)}
\Ofa{\p3k}(-\delta)\to\Ofa F(-\delta)\to0.\tag 1.5
$$
Applying $\Hom{{\p3k}}\big({\Cal F},\cdot\big)$ to sequence
(1.5) one gets
$$
0\to\Hom{F}\big({\Cal F},\Ofa
F(-\delta)\big)\to\Ext1{\p3k}\big({\Cal
F},\Ofa{\p3k}(-d-\delta)\big)\mapright{f}
\Ext1{\p3k}\big({\Cal F},\Ofa{\p3k}(-\delta)\big),
$$
and since the multiplication by $f$ is
zero on
$\Ext1{\p3k}\big({\Cal F},\Ofa{\p3k}(-d-\delta)\big)$ we
obtain an isomorphism $s_0\colon{\Cal F}\mapright\sim
\Hom{F}\big({\Cal F},\Ofa
F(-\delta)\big)\to\Ext1{\p3k}\big({\Cal
F},\Ofa{\p3k}(-d-\delta)\big)$.

\proclaim{Proposition 1.6}
Let $s_0$ be as above. If there is a commutative diagram
$$
\matrix
0&\hskip-3truemm\to&\hskip-3truemm{\Cal
K}&\hskip-3truemm\mapright{d}&\hskip-3truemm{\Cal
E}&\hskip-3truemm\mapright{\gamma}&\hskip-3truemm{\Cal
F}&\hskip-3truemm\to&\hskip-3truemm0\\
&\hskip-3truemm&\hskip-3truemm\mapdown{s_2}
&\hskip-3truemm&\hskip-3truemm\mapdown{s_1}
&\hskip-3truemm&\hskip-3truemm\mapdown{s_0}
&\hskip-3truemm&\hskip-3truemm\\
0&\hskip-3truemm\to&\hskip-3truemm \check{\Cal E}(-d-\delta)&
\hskip-3truemm\mapright{\check{d}}& \hskip-3truemm\check{\Cal
K}(-d-\delta)& \hskip-3truemm\mapright{\overline{\gamma}}
&\hskip-3truemm\Ext1{\p3k}\big({\Cal
F},\Ofa{\p3k}(-d-\delta)\big)
&\hskip-3truemm\to&\hskip-3truemm0.
\endmatrix\tag 1.6.1
$$
then the maps $s_i$ are isomorphisms.
\endproclaim
\demo{Proof}
By construction $s_0$ is an isomorphism thus it suffices to
prove the bijectivity of $s_1$. This will follow from lemma
2.12 of [Wa] if we show that conditions (i), (ii) and (iii)
are verified for the composition
$$
{\Cal E}\mapright{s_1}\check{\Cal
K}(-d-\delta)
\mapright{\overline{\gamma}}\Ext1{\p3k}\big({\Cal
F},\Ofa{\p3k}(-d-\delta)\big).
$$
Condition (iii) that $\Cal E$ and $\check{\Cal K}(-d-\delta)$
have the same rank was already shown in proposition 1.4.

Since $\overline{\gamma}\circ s_1$ coincides with the map
$\gamma\colon {\Cal E}\twoheadrightarrow{\Cal F}$ condition
(ii) holds by the very definition of $\Cal E$.

We have to verify (i) namely that
$$
H^i_*\big(s_1\big)\colon H^i_*\big(\p3k,{\Cal E}\big)\to
H^i_*\big(\p3k,\check{\Cal K}(-d-\delta)\big)
$$
are isomorphisms for $i=1,2$. Diagram (1.6.1) yields the
equality
$$
H^i_*\big(\overline\gamma\big)\circ H^i_*\big(s_1\big)=
H^i_*\big(s_0\big)\circ H^i_*\big(\gamma\big).
$$
Note that $H^i_*\big(s_0\big)$ are isomorphisms and the
maps $H^i_*\big(\gamma\big)$ are injective hence the same is
true for the maps
$H^i_*\big(s_1\big)$.

Since $\Cal K$ and $\Cal E$ are locally free then
both $H^i_*\big(\p3k,\check{\Cal K}(-d-\delta)\big)$ and
$H^i_*\big(\p3k,{\Cal E}\big)$ have finite length thus we
have only to prove that their  lengths coincide.

We begin with $i=2$. Here both modules are $0$:
$$
h^2\big(\p3k,\check{\Cal K}(t)\big)=h^1\big(\p3k,{\Cal
K}(-4-t)\big)=0
$$
by iv) while $h^2\big(\p3k,{\Cal E}(t)\big)=0$ (cf. iii)).

Let now $i=1$. One has
$$
\align
h^1\big(\p3k,\check{\Cal K}(t)\big)&
=\dim\left(H^2\big(\p3k,{\Cal
K}(-4-t)\big)\check{}\right)=\dim\left(H^1\big(F,{\Cal
F}(-4-t)\big)/W_{-4-t}\right)\cong\\ &=\cases
0&\text{ if $-4-t<(d-4+\delta)/2$},\\
\dim(U)&\text{ if $-4-t=(d-4+\delta)/2$},\\
h^1\big(F,{\Cal F}(-4-t)\big)&\text{ if
$-4-t>(d-4+\delta)/2$}.
\endcases
\endalign
$$
By (1.2) $h^1\big(\p3k,{\Cal F}(-4-t)\big)
=h^1\big(\p3k,{\Cal F}(t+d+\delta)\big)$ and so the desired
conclusion follows.
\qed
\enddemo

\head
2. Proof of theorem 0.3
\endhead

This section is devoted to the proof of the following

\proclaim{Claim 2.1}
It is possible to construct diagram (1.6.1) in
such a way that $s_2$ is the transpose of $s_1$.
\endproclaim

Assuming 2.1 we have the

\demo{Proof of theorem 0.3}
Just set $\varphi:=s_1^{-1}\circ\check d$ which is
obviously symmetric.
\qed
\enddemo

\demo{Proof of claim 2.1}
Our first step is to extend the natural map
$\eta\colon\Sim^2{\Cal F}\to\Ofa F(-\delta)$, induced by the
symmetric map $\sigma$, to a chain map $\phi\colon
\Sim^2{\Cal R}^*\to{\Cal C}^*$ (see sequences (1.3) and (1.5))
$$
\matrix
0&\hskip-3truemm\to&\hskip-3truemm\Lambda^2{\Cal
K}&\hskip-3truemm\mapright{\delta_1}&\hskip-3truemm{\Cal
K}\otimes{\Cal
E}&\hskip-3truemm\mapright{\delta_0}&
\hskip-3truemm\Sim^2{\Cal
E}&\hskip-3truemm\to&\hskip-3truemm\Sim^2{\Cal
F}&\hskip-3truemm\to&\hskip-3truemm0\\
&\hskip-3truemm&\hskip-3truemm\mapdown{\phi_2}
&\hskip-3truemm&\hskip-3truemm\mapdown{\phi_1}
&\hskip-3truemm&\hskip-3truemm\mapdown{\phi_0}
&\hskip-3truemm&\hskip-3truemm\downarrow&\hskip-3truemm&
\hskip-3truemm\\
&\hskip-3truemm&\hskip-3truemm0&\hskip-3truemm\to
&\hskip-3truemm\Ofa{\p3k}(-d-\delta)&
\hskip-3truemm\mapright{\lambda(f)}
&\hskip-3truemm\Ofa{\p3k}(-\delta)&\hskip-3truemm\to&
\hskip-3truemm\Ofa
F(-\delta)&\hskip-3truemm\to&\hskip-3truemm0.
\endmatrix
$$
Assume that $\phi$ does exist. Then we get a map
$s_1\colon{\Cal E}\to\check{\Cal K}(-d-\delta)$. It is
obtained from $\phi_1$ through the natural isomorphism
$$
\Hom{{\p3k}}\big({\Cal K}\otimes{\Cal
E},\Ofa{\p3k}(m)\big)\cong\Hom{{\p3k}}\big({\Cal
E},\check{\Cal K}(m)\big).
$$
Let $s_2$ be the transpose of $s_1$.

We claim that $s_1\circ d=\check
d\circ s_2$ i.e. that the diagram (1.6.1) above actually
commutes. It suffices to verify this equality at every point
$x\in\p3k$. Choose
$\alpha,\beta\in{\Cal K}_x$. Since $\langle
s_1\circ d(\alpha),\beta\rangle=
\phi_1(d(\alpha)\otimes\beta)$ and
$\langle
\check d\circ s_2(\alpha),\beta\rangle=\langle
\alpha,s_1\circ d(\beta)\rangle= \phi_1(\alpha\otimes
d(\beta))$, our claim follows
$$
0=\phi_2(\alpha\wedge\beta)=
\phi_1\circ\delta_1(\alpha\wedge\beta)=
\phi_1(\alpha\otimes d(\beta)-\beta\otimes d(\alpha)).
$$
There remains only to prove the following proposition.

\proclaim{Proposition 2.2}
$\phi$ exists.
\endproclaim
\demo{Proof}
In order to have $\phi$ it suffices to define $\phi_0$.
Indeed the image of $\phi_0\circ\delta_0$ is contained in the
kernel of
$\Ofa{\p3k}\to\Ofa F$ which coincides with the ideal
$f\Ofa{\p3k}$.

Then we simply set $\phi_1:={1\over f}\phi_0\circ\delta_0$.

We want to lift the composition
$\psi:=\eta\circ\Sim^2\gamma\colon
\Sim^2{\Cal E}\to\Sim^2{\Cal F}\to\Ofa F(-\delta)$ to a
map $\phi_0\colon\Sim^2{\Cal E}\to\Ofa{\p3k}(-\delta)$.

{}From the exact sequence (1.5) we obtain the exact sequence
$$
\align
0&\to\hom_{\Ofa{\p3k}}\big(\Sim^2{\Cal
E},\Ofa{\p3k}(-d-\delta)\big)\to
\hom_{\Ofa{\p3k}}\big(\Sim^2{\Cal
E},\Ofa{\p3k}(-\delta)\big)\to\\
&\to\hom_{\Ofa{\p3k}}\big(\Sim^2{\Cal E},\Ofa
F(-\delta)\big)
\mapright\partial\ext^1_{\Ofa{\p3k}}\big(\Sim^2{\Cal
E},\Ofa{\p3k}(-d-\delta)\big)\cong\\ &\cong
H^1\big(\p3k,\Sim^2\check{\Cal E}(-d-\delta)\big)\cong
H^2\big(\p3k,\Sim^2{\Cal E}(d-4+\delta)\big)\check{}.
\endalign
$$
We conclude that $\psi$ is liftable if and only if
$\partial(\psi)\in H^2\big(\p3k,\Sim^2{\Cal
E}(d-4+\delta)\big)\check{}$ is the zero map. First of all
notice that, interchanging the roles of $\Lambda^2$ and
$\Sim^2$ in section 4 of [Wa] and imitating word by word the
proofs of lemmas 4.1, 4.2 and corollary 4.3 of [Wa] we
easily obtain
$$
H^2\big(\p3k,\Sim^2{\Cal E}(d-4+\delta)\big)\cong\cases
0&\text{ if $d+\delta$ is odd,}\\
\Lambda^2 H^1\big(\p3k,{\Cal E}((d-4+\delta)/2)\big)& \text{
if
$d+\delta$ is even.}
\endcases
$$

$\partial(\psi)\in H^2\big(\p3k,\Sim^2{\Cal
E}(d-4+\delta)\big)\check{}$ is identified with the map
$$
\big\lrcorner(\cdot\smallsmile\partial(\psi))\colon
H^2\big(\p3k,\Sim^2{\Cal E}(d-4+\delta)\big)\to
H^3\big(\p3k,\Ofa{\p3k}(-4)\big),
$$
where $\smallsmile\colon H^2\big(\p3k,\Sim^2{\Cal
E}(d-4+\delta)\big)\times H^1\big(\p3k,\Sim^2\check{\Cal
E}(-d-\delta)\big)\to H^3\big(\p3k,\Sim^2{\Cal
E}\otimes\Sim^2\check{\Cal E}(-4)\big)$ is the cup--product
and $\big
\lrcorner\colon H^3\big(\p3k,\Sim^2{\Cal
E}\otimes\Sim^2\check{\Cal E}(-4)\big)\to
H^3\big(\p3k,\Ofa{\p3k}(-4)\big)$ is the natural contraction.

Thus $\partial(\psi)=0$ if and only if
$$
\lrcorner(\alpha\smallsmile\beta\smallsmile\partial(\psi))=0
\tag2.2.1
$$
in $H^3\big(\p3k,\Ofa{\p3k}(-4)\big)$ for each
$\alpha,\beta\in H^1\big(\p3k,{\Cal E}((d-4+\delta)/2)\big)$.

We have a commutative diagram
$$
\matrix
\Lambda^2 H^1\big(\p3k,{\Cal
E}((d-4+\delta)/2)\big)&\hskip-3truemm
\mapright{\Lambda^2H^1\big(\gamma\big)}&\hskip-3truemm
\Lambda^2H^1\big(F,{\Cal
F}((d-4+\delta)/2)\big)&\hskip-3truemm&\hskip-3truemm\\
\Vert&\hskip-3truemm&\hskip-3truemm\mapdown\Phi
&\hskip-3truemm&\hskip-3truemm\\
H^2\big(\p3k,\Sim^2{\Cal
E}(d-4+\delta)\big)&\hskip-3truemm
\mapright{H^2\big(\psi\big)}&\hskip-3truemm H^2\big(F,\Ofa
F(d-4)\big)&\hskip-3truemm\mapright{\partial'}
&\hskip-3truemmH^3\big(\p3k,\Ofa{\p3k}(-4)\big)
\endmatrix
$$

We claim that
\proclaim{Lemma 2.3}
$\big\lrcorner(\cdot\smallsmile\partial(\psi))
=-\partial'\circ H^2\big(\psi\big)$.
\endproclaim
The above assertion implies $\partial(\psi)=0$ since, by
formula 2.2.1, $\partial'\circ
H^2\big(\psi\big)(\alpha\smallsmile\beta)=0$ because
$\alpha,\beta\in H^1\big(\p3k,{\Cal
E}((d-4+\delta)/2)\big)=U$ which was chosen to be isotropic
with respect to $\Phi$.

\demo{Proof of lemma 2.3} Let ${\Cal U}:=\lbrace
U_i\rbrace_{i=0,\dots,3}$ be the standard open covering of
$\p3k$. For each $i$ we fix a lifting
$\widehat{\psi}_i\colon\Sim^2{\Cal E}_{\vert
U_i}\to\Ofa{U_i}(-\delta)$ of $\psi_{\vert U_i}$. Notice that
$\widehat\psi_i-\widehat\psi_j$ maps to
$f\Ofa{\p3k}(-d-\delta)\subseteq\Ofa{\p3k}(-\delta)$. On the
other hand
$$
\partial(\psi)\in H^2\big(\p3k,\Sim^2{\Cal
E}(d-4+\delta)\big)\check{\ }\cong
H^1\big(\p3k,\Sim^2\check{\Cal E}(-d-\delta)\big)
$$
represents the obstruction to lifting $\psi$ to $\phi_0$
whence
$$
\partial(\psi)={1\over f}(\widehat\psi_i-\widehat\psi_j)\in
H^1\big({\Cal U},\Sim^2\check{\Cal E}(-d-\delta)\big).
$$

We now compute explicitly $\partial'\circ
H^2\big(\psi\big)$ and
$\big\lrcorner(\cdot\smallsmile\partial(\psi))$ using the
fact that each element inside
$H^2\big(\p3k,\Sim^2{\Cal E}(d-4)\big)$ can be written as a
sum
$\alpha\smallsmile\beta$ where
$\alpha,\beta\in H^1\big({\Cal U},{\Cal E}((d-4)/2)\big)$.
$\alpha\smallsmile\beta\smallsmile\partial(\psi)$ is
represented by the cocycle
$$
\big(\alpha\smallsmile\beta\smallsmile\partial(\psi)
\big)_{i_0,i_1,i_2,i_3}= {1\over
f}(\alpha_{i_0,i_1}\beta_{i_1,i_2})\otimes
(\widehat\psi_{i_2}-\widehat\psi_{i_3}),
$$
hence
$$
\big(\big\lrcorner
(\alpha\smallsmile\beta\smallsmile\partial(
\psi))\big)_{i_0,i_1,i_2,i_3} ={1\over
f}\big(\widehat\psi_{i_2}(\alpha_{i_0,i_1}\beta_{i_1,i_2})-
\widehat\psi_{i_3}(\alpha_{i_0,i_1}\beta_{i_1,i_2})\big).
$$
On the other hand
$H^2\big(\psi\big)(\alpha\smallsmile\beta)_{i,j,h}=
\psi(\alpha_{i,j}\beta_{j,h})\psi$ whence
$$
\align
\partial'\circ
H^2\big(\psi\big)(
\alpha\smallsmile\beta)_{i_0,i_1,i_2,i_3}&= {1\over f}
\big(\widehat\psi_{j_0}(\alpha_{i_1,i_2}\beta_{i_2,i_3})
-\widehat\psi_{j_1}(\alpha_{i_0,i_2}\beta_{i_2,i_3})+\\
&+\widehat\psi_{j_2}(\alpha_{i_0,i_1}\beta_{i_1,i_3})
-\widehat\psi_{j_3}(\alpha_{i_0,i_1}\beta_{i_1,i_2})\big),
\endalign
$$
where $j_h\in\lbrace i_0,i_1,i_2,i_3\rbrace\setminus\lbrace
i_h\rbrace$. In particular choosing  $j_0=j_1=j_2=i_3$ and
$j_3=i_2$ and using that
$\alpha_{i_1,i_2}-\alpha_{i_0,i_2}+\alpha_{i_0,i_1}=0$,
$\beta_{i_1,i_2}-\beta_{i_1,i_3}+\beta_{i_2,i_3}=0$ we get
$$
\align
\partial'\circ
H^2\big(\psi\big)(
\alpha\smallsmile\beta)_{i_0,i_1,i_2,i_3}&= {1\over f}
\big(\widehat\psi_{i_3}(-\alpha_{i_0,i_1}(\beta_{i_2,i_3}-
\beta_{i_1,i_3}))
-\widehat\psi_{i_2}(\alpha_{i_0,i_1}\beta_{i_1,i_2})\big)=\\
&={1\over f}
\big(\widehat\psi_{i_3}(\alpha_{i_0,i_1}\beta_{i_1,i_2})
-\widehat\psi_{i_2}(\alpha_{i_0,i_1}\beta_{i_1,i_2})\big).
\endalign
$$
Then the proof is complete.
\qed
\qed
\qed
\enddemo
\enddemo
\enddemo

\remark{Remark 2.2}
Theorem 0.3 holds without the assumption $F\subseteq\p3k$.
It suffices to consider any hypersurface $F\subseteq\p nk$
endowed with a
$\delta/2$--quadratic sheaf $\Cal F$.
\endremark

\head
3. Examples and applications
\endhead

In this section we shall consider the case of a nodal
surface of small degree
$d=4,5$ and we shall see how our main theorem can be used to
classify $\delta/2$--even sets of nodes, going beyond
[Ct], section 3 and [Bea].

Moreover from now on the ground field $k$ will be equal to
the field
$\Bbb C$ of complex numbers.

Following the notations used in the introduction there exists
an invertible $\Ofa {\widetilde{F}}$--sheaf $\Cal L$ such
that
$\pi_*{\Cal L}\cong\Cal F$,
$\widetilde{p}_*\Ofa{S}\cong\Ofa{\widetilde{F}}\oplus{\Cal
L}$. Since
$R^1\widetilde{\pi}_*\Ofa{\widetilde{S}}=
R^1\pi_*\Ofa{\widetilde{F}}=0$, diagram (0.1) and the
spectral sequences of the composite functors
$R^p\pi_*R^q\widetilde{p}_*$ and
$R^pp_*R^q\widetilde{\pi}_*$ yields
$R^1\pi_*{\Cal L}=0$ thus
$h^i\big(\widetilde{F},{\Cal L}\otimes\pi^*\Ofa
F(m)\big)=h^i\big(F,{\Cal F}(m)\big)$.

Note that
$$
h^1\big(F,{\Cal F}(-m)\big)=h^1\big(F,{\Cal
F}(m+d-4+\delta)\big)=0,\quad m>0\tag3.1
$$
(theorem 1 of [Kw] applied to
$\widetilde{p}^*\pi^*\Ofa{\widetilde{S}}(1)$ and formula
(1.2)). Moreover by [Ct],
theorem 2.19, we get that $\Delta$ is symmetric if and only
if
$h^1\big(F,{\Cal F}(m)\big)=0$ for $0\le m\le(d-4)/2$.

We have  long exact sequence
$$
\aligned
0&\to H^0\big({F},{\Cal F}\big)\to H^0\big({F},{\Cal
F}(1)\big)\to H^0\big({H},{\Cal F}(1)_{\vert {H}}\big)\to
H^1\big({F},{\Cal F}\big)\mapright{\lambda(h)}\\
&\to H^1\big({F},{\Cal F}(1)\big)\to
H^1\big({H},{\Cal F}(1)_{\vert {H}}\big)\to
H^2\big({F},{\Cal F}\big)\to H^2\big({F},{\Cal
F}(1)\big)\to0.
\endaligned\tag3.2
$$
associated to every $h\in H^0\big(\p3k,\Ofa{\p3k}(1)\big)$
(defining a plane section $H$ of $F$). Notice that in any
case $h^0\big(F,{\Cal F}\big)=0$ since $S$ is connected.

Moreover, by (3.1), then $H^1\big(F,{\Cal F}(-1)\big)=0$,
thus
$H^1\big(H,{\Cal F}_{\vert H}\big)=0$ for each $H$.

Since $h$ does not divide $f:=\det(\varphi)$, the functor
$\cdot\otimes_{\Ofa{F}}\Ofa H$ is exact. It follows by the
exact sequence 0.3.1 that ${\Cal F}_{\vert H}$ is again a
$\delta/2$--quadratic sheaf.
Mutatis mutandis, theorem 2.16 and proposition 2.28 of [Ct]
apply, therefore
${\Cal F}_{\vert H}$ has a free resolution which is exact on
global sections
$$
0\to\bigoplus_{j=1}^h\Ofa
H(-\ell_j) \mapright\alpha\bigoplus_{i=1}^h\Ofa
H(-r_i)\to{\Cal F}_{\vert H}\to0\tag 3.3
$$
where $\alpha:=\left(\alpha_{i,j}\right)_{i,j=1,\dots,h}$ is
a symmetric matrix of homogeneous polynomials $\alpha_{i,j}$
of degrees
$(d_i+d_j)/2$, where the $d_i$'s are in not--decreasing
order,
$d_i\equiv d_j,\ d\equiv\delta+d_i\pmod2$,
$\ell_j=(d+\delta+d_j)/2$ and
$r_i=(d+\delta-d_i)/2$.

As in [Ct] we see that:
\roster
\item"{i)}" $d_i+d_{h+1-i}>0$ since $\det(\alpha)$;
\item"{ii)}" $d_i+d_{h-i}>0$ if $\det(\alpha)$ is square
free;
\item"{iii)}" $r_i>0$ since $h^0\big(H,{\Cal F}_{\vert
H}\big)=0$, i.e.
$d_i\le d+\delta-2$;
\item"{iv)}" $d_i+d_{h-1-i}>0$ if $H=\lbrace
\det(\alpha)=0\rbrace$ is smooth.
\endroster

Notice finally that
$$
d=\sum_{i=1}^hd_i=\sum_{i\le h/2}(d_i+d_{h+1-i})+d_{(h+1)/2}.
$$
Here $d_\lambda=0$ if $\lambda$ is not an integer. We then
get the following cases for the
$h$--tuple $(d_1,\dots,d_h)$:
\roster
\item"{$d=4$, $\delta=0$}" $(2,2)$, $(0,2,2)$, $(0,0,2,2)$;
\item"{$d=4$, $\delta=1$}" $(1,3)$, $(1,1,1,1)$,
$(-1,1,1,3)$, $(-1,-1,3,3)$;
\item"{$d=5$, $\delta=0$}" $(1,1,3)$, $(-1,3,3)$, $(-1,1,5)$,
$(1,1,1,1,1)$, $(-1,1,1,1,3)$, $(-1,-1,1,3,3)$.
\endroster
On the other hand if we assume that $H$ is smooth most of
the above possibilities disappear and we are then only left
with:
\roster
\item"{$d=4$, $\delta=0$}" $(2,2)$;
\item"{$d=4$, $\delta=1$}" $(1,3)$, $(1,1,1,1)$;
\item"{$d=5$, $\delta=0$}" $(1,1,3)$, $(1,1,1,1,1)$.
\endroster

\subhead
3.4. The case $d=4$, $\delta=0$
\endsubhead
For all $H$ we have $h^1\big({H},{\Cal
F}(1)_{\vert{H}}\big)=h^0\big({H},{\Cal
F}_{\vert{H}}\big)=0$ then $h^0\big({H},{\Cal
F}(1)_{\vert{H}}\big)=2$ (Riemann--Roch) hence formula (3.1)
yields
$h^1\big({F},{\Cal F}(1)\big)=h^1\big({F},{\Cal
F}(-1)\big)=0$. Moreover
$h^1\big({F},{\Cal F}\big)$ is dual to
itself whence it has even dimension. Recall that
$$
\chi\big({\Cal F}\big)=(8-t)/4
$$
(see [Ct], proposition 2.11).

If $h^1\big({F},{\Cal F}\big)=0$ then $\Delta$ is symmetric
of type $(2,2)$,
$h^0\big({F},{\Cal F}(1)\big)=2$ and $t=8$.

If $h^1\big({F},{\Cal F}\big)=2$, $h^0\big({F},{\Cal
F}(1)\big)=0$, $t=16$. It follows from
condition ii) of section 1 that $h^1\big(\p3k,{\Cal
E}\big)=1$ and
$h^1\big(\p3k,{\Cal E}(m)\big)=0$ for $m\ne0$.

The Horrocks correspondence shows that $\Cal E$ is stably
equivalent to
$\Omega^1_{\p3k}$ (stably equivalent means that adding a
direct sum of invertible sheaves we get isomorphic sheaves).
Therefore, by the construction of $\Cal E$, since
$h^0\big(\p3k,\Omega^1_{\p3k}\big)=
h^0\big(\p3k,\Omega^1_{\p3k}(1)\big)=0$, we may assume that
$h^0\big(\p3k,{\Cal E}(1)\big)=0$. Since we have
$h^i\big(F,{\Cal F}(1)\big)=0$ for $i=0,1$, $h^0\big(F,{\Cal
F}(2)\big)=h^0\big(H,{\Cal F}(2)_{\vert
H}\big)=h^0\big(\p3k,\Omega^1_{\p3k}(2)\big)=6$. Recall
(section 1 ii)) that there is an epimorphism
$H^0\big(\p3k,{\Cal E}(2)\big)\twoheadrightarrow
H^0\big(\p3k,{\Cal F}(2)\big)$. Let $6-r$ be the rank of the
induced map
$H^0\big(\p3k,\Omega^1_{\p3k}(2)\big)\to H^0\big(\p3k,{\Cal
F}(2)\big)$. Then
$$
{\Cal E}\cong\Omega_\p3k^1\oplus\Ofa{\p3k}(-2)^{\oplus
r}\oplus\Ofa{\p3k}(-3)^{\oplus r_3}\oplus\dots.
$$
We claim that $r_3=r_4=\dots=0$.

Indeed, by Beilinson's theorem (see [Bei]) applied to ${\Cal
E}(1)$, since $H^3\big(\p3k,{\Cal E}(-1)\big)\check{\
}=H^0\big(\p3k,\check{\Cal E}(-3)\big)$ which, by theorem
0.3, injects into $H^0\big(\p3k,{\Cal E}(1)\big)=0$, we get
$$
{\Cal E}(1)\cong\left(H^1\big(\p3k,{\Cal
E}\big)\otimes\Omega_\p3k^1(1)\right)
\oplus\left(H^3\big(\p3k,{\Cal
E}(-2)\big)\otimes\Ofa{\p3k}(-1)\right),
$$
thus ${\Cal E}\cong\Omega_\p3k^1
\oplus\Ofa{\p3k}(-2)^{\oplus r}$ and
$\check{\Cal E}(-4)\cong\Omega_\p3k^2
\oplus\Ofa{\p3k}(-2)^{\oplus r}$. It follows that
$$
\varphi=\pmatrix
\varphi_{1,1}&\varphi_{1,2}\\
\varphi_{1,2}&\varphi_{2,2}
\endpmatrix
$$
where $\varphi_{2,2}$ is a constant $r\times r$ symmetric
matrix which must be zero by the minimality of $\Cal E$.

If $r\ge4$ then $\det(\varphi)=0$ which is absurd.

If $r=3$ then $\rk(\varphi_{1,2})\le2$ on a surface $Y$ of
degree
$c_1\big(\Omega_\p3k^1(2)\big)=2$, then $F=2Y$, again an
absurd.

If $r=2$ there is a curve $\Gamma$ of degree
$c_2\big(\Omega_\p3k^2(2)\big)=2$ where
$\rk(\varphi_{1,2})\le1$, whence
$\Gamma$ is a double curve for $F$, again an absurd.

If $r=1$ then $\varphi_{1,2}$ is a section of
$H^0\big(\p3k,\Omega_\p3k^1(2)\big)\cong\Lambda^2
H^0\big(\p3k,\Ofa{\p3k}(1)\big)$.

If the rank of this alternating map
is $2$ then $\varphi_{1,2}$ should vanish on a line $\Gamma$
(since $\varphi_{1,2}=x_0dx_1-x_1dx_0$ for suitable
coordinates) which is a double line for $F$, absurd.

We conclude that $\varphi_{1,2}$
corresponds to a non--degenerate alternating form. With a
proper choice of the coordinates we can assume that
$\varphi_{1,2}$ corresponds to
$x_0dx_1-x_1dx_0+x_2dx_3-x_3dx_2$. Then $\Ofa{\p3k}(-2)$ is
identified to a subbundle of $\Omega_\p3k^1$, and dually
$\Omega_\p3k^2\twoheadrightarrow\Ofa{\p3k}(-2)$ is
surjective. As in [Ba] we define the null--correlation bundle
${\Cal V}_0$ as
${\Cal V}_0:=\Omega_\p3k^1(2)/\im(\varphi_{1,2}(2))$ and we
obtain a self--dual resolution
$$
0\to\check{\Cal V}_0\to{\Cal V}_0\to{\Cal F}(2)\to0.
$$

Finally If $r=0$ we get
$$
0\to\big(\Omega_\p3k^1(2)\big)\check{}
\to\Omega_\p3k^1(2)\to{\Cal F}(2)\to0.
$$

We have therefore shown that for quartic surfaces all even
sets of nodes are either symmetric or $t=16$ and we have
exactly the two cases described in [Ba1].

\subhead
3.5. The case $d=4$, $\delta=1$
\endsubhead
As in the case $\delta=0$ we have again $h^1\big({F},{\Cal
F}(1)\big)=0$. We get by (1.2) that $h^1\big({F},{\Cal
F}\big)=h^1\big({F},{\Cal
F}(1)\big)=0$. It then follows
that $\Delta$ is symmetric of type either
$(1,1,1,1)$ or $(1,3)$.

\subhead
3.6. The case $d=5$
\endsubhead
In this case $\delta=0$ as we already noticed. Sequence
(3.2) is self-dual and
$$
\lambda\colon
H^0\big(\p3k,\Ofa{\p3k}(1)\big)
\to\hom_k\big(H^1\big({F},{\Cal F}(1)\big),H^1\big({F},{\Cal
F}(1)\big)\check{}\big)
$$
maps to the subspace of alternating bilinear forms.

If ${H}$ is smooth ${\Cal F}_{\vert H}$ is of type either
$(1,1,1,1,1)$ or $(1,1,3)$ thus $h^0\big({H},{\Cal
F}(1)_{\vert{H}}\big)= h^1\big({H},{\Cal
F}(1)_{\vert{H}}\big)\le 1$ by (3.3).

The map $\lambda(h)$ is an isomorphism if either
$h^0\big({H},{\Cal F}(1)_{\vert{H}}\big)=0$ or
$h^0\big({F},{\Cal F}(1)\big)=1$. Infact in both the cases
then $h^0\big(F,{\Cal F}(1)\big)=h^0\big(H,{\Cal F}(1)_{\vert
H}\big)$ and the assertion follows easily from the
self--duality of sequence (3.2).

If for a general $H$ the map $\lambda(h)$ is an isomorphism
then
$h^1\big(F,{\Cal F}\big)$ is even and if it is not zero the
pfaffian of
$\lambda(h)$ defines a surface $B\subseteq\check\p3k$ of
degree
$h^1\big({F},{\Cal F}\big)/2$ which is
contained in the dual surface $\check F$ of $F$. Since $F$
is nodal then
$\check F$ is of general type and, by biduality, one has
$\deg(\check F)\ge5$, whence $h^1\big({F},{\Cal
F}\big)\ge10$. We conclude that
$$
t/4-5=-\chi({\Cal F})\ge 10
$$
hence $t\ge60$ which implies $\deg(\check{F})\le
d(d-1)^2-2t\le-40$, an absurd. Thus the only possibility is
$h^1\big({F},{\Cal F}\big)=0$ and we get that
$\Delta$ is symmetric.

Finally if, for $H$ smooth, $h^0\big({H},{\Cal
F}(1)_{\vert{H}}\big)=1$ and $h^0\big({F},{\Cal
F}(1)\big)=0$ then $\dim(\ker(\lambda(h)))=1$, thus
$H^1\big({F},{\Cal F}\big)\ne0$ is odd ($\lambda(h)$ is
alternating).

In any case $h^1\big(F,{\Cal F}(2)\big)=0$ (again
(3.1)) hence we have the exact sequence
$$
0\to H^0\big(F,{\Cal F}(2)\big)\to H^0\big(H,{\Cal
F}(2)_{\vert H}\big)\to H^1\big(F,{\Cal F}(1)\big)\to0.
$$
In particular $h^0\big(F,{\Cal F}(2)_{\vert H}\big)$ does
not depend on
$H$. If $H$ is general, either of type $(1,1,1,1,1)$ or
$(1,1,3)$, sequence (3.3) implies $h^0\big(H,{\Cal
F}(2)_{\vert H}\big)=5$. On the other hand if
$h^1\big({F},{\Cal F}\big)\ge3$ there exists a plane section
$H$ such that $\dim(\ker(\lambda(h)))\ge3$, whence of type
$(-1,1,5)$. Thus looking at sequence (3.3) we get
$h^0\big(H,{\Cal F}(2)_{\vert H}\big)=6$, a contradiction.

We have restricted ourselves to the case $h^1\big({F},{\Cal
F}\big)=1$.

In this case
$h^1\big({F},{\Cal F}(-1)\big)=0$ (formula (3.1)) and
$h^0\big(F,{\Cal F}(-1)\big)=0$ thus
$h^0\big(F,{\Cal F}(2)\big)=h^2\big(F,{\Cal F}(-1)\big)=4$.
Beilinson's theorem then yields a
sequence of the form
$$
0\to\Ofa{\p3k}(-1)^{\oplus4}\oplus\Omega_{\p3k}^2(2)
\mapright\varphi{\Cal
O}_{\p3k}^{\oplus4} \oplus\Omega_{\p3k}^1(1)\to{\Cal
F}(2)\to0,
$$
where
$$
\varphi=\pmatrix
\varphi_{1,1}&\varphi_{1,2}\\
\varphi_{2,1}&\varphi_{2,2}
\endpmatrix.
$$
Moreover Beilinson's spectral sequence and the fact that
$\lambda(h)=0$ for each $h$ implies that
$\varphi_{2,2}\colon\Omega_{\p3k}^2(2)\to\Omega_{\p3k}^1(1)$
is zero. In particular we see that ${\Cal
O}_{\p3k}^{\oplus4}\oplus\Omega_{\p3k}^1(1)$ has the
required properties for $\Cal E$, thus we may assume that
$$
{\Cal E}(2)\cong{\Cal
O}_{\p3k}^{\oplus4}\oplus\Omega_{\p3k}^1(1)
$$
and that $\varphi$ is symmetric.

But then the determinantal quintic $F$ should be singular
along the set
$D\subseteq\p3k$ of points where $\rk(\varphi_{1,2})\le2$
which has dimension at least $1$ (generically $D$ is a pair
of skew lines). Infact
$\det(\varphi)$ belongs to the square of the sheaf of ideals
of $D$.

\enddocument